\documentclass[10pt,reqno]{amsart}
\usepackage{amsfonts,amsmath,amssymb,amsthm}
\usepackage{latexsym}
\usepackage{eucal}
\usepackage[ansinew]{inputenc}
\usepackage[american]{babel}
\usepackage{amsfonts}
\usepackage{amssymb}
\usepackage{graphicx}

\newtheorem{teo}{Theorem}[section]

\newtheorem{defi}{Definition}[section]
\newtheorem{obs}{Remark}[section]

\begin{document}

\title[Stability for the Kawahara-KdV Type Equations]
{A note on the stability for Kawahara-KdV type equations.}

\author[F. Natali]{F\'abio Natali $^1$}

\email{fmnatali@hotmail.com, fmanatali@uem.br}

\thanks{{\it Date}: June, 2009.}

\maketitle

{\scriptsize \centerline{$^1$Department of Mathematics, Universidade
Estadual de Maringá}
 \centerline{Avenida Colombo, 5790, CEP 87020-900, Maringá, PR, Brazil.}}
%%%%%%%%%%%%%%%%%%%%%%%%%%%%%%%%%%%%%%%%%%%%

\begin{abstract} In this paper we establish the nonlinear stability of
solitary traveling-wave solutions for the Kawahara-KdV equation
$$u_t+uu_x+u_{xxx}-\gamma_1 u_{xxxxx}=0,$$
and the modified Kawahara-KdV equation
$$u_t+3u^2u_x+u_{xxx}-\gamma_2 u_{xxxxx}=0,$$
where $\gamma_i\in\mathbb{R}$ is a positive number when $i=1,2$. The
main approach used to determine the stability of solitary
traveling-waves will be the theory developed by Albert in
\cite{albert1}.

\end{abstract}
\section{Introduction.}
This work presents the existence of a smooth branch of solitary
traveling wave solutions as well as the orbital stability related to
Kawahara-Korteweg-de Vries and modified Kawahara-Korteweg-de Vries
equations (Kawahara and modified Kawahara equations respectively,
henceforth),
\begin{equation}
u_t+uu_x+u_{xxx}-\gamma_1 u_{xxxxx}=0,
\label{equakawa}\end{equation} and
\begin{equation}
u_t+3u^2u_x+u_{xxx}-\gamma_2 u_{xxxxx}=0,
\label{equamkawa}\end{equation} where $\gamma_i>0$ when $i=1,2$ and
$u:=u(x,t)$ is a real function. Here, we consider $x\in\mathbb{R}$
and $t\in\mathbb{R}$. These equations model the propagation on
nonlinear water-waves in the long-wavelength as in the case KdV's
equations. Roughly speaking, such a model-scenario is expected
because, if $u$ be a smooth solution of $(\ref{equakawa})$ and
$(\ref{equamkawa})$, then for $\gamma_i\rightarrow0$ uniformly,
$i=1,2$ we obtain that $u$ is a solution of the Korteweg-de Vries
and modified Korteweg-de Vries equations,
\begin{equation}
u_t+uu_x+u_{xxx}=0, \label{kdv}\end{equation}

\begin{equation}
u_t+3u^2u_x+u_{xxx}=0, \label{mkdv}\end{equation} respectively, in a
convenient sense. Results of orbital stability for equations
$(\ref{kdv})$ and $(\ref{mkdv})$ has been studied by many
researchers in the case of solitary waves, for example see
\cite{albert1}, \cite{albert2}, \cite{albert3}, \cite{benjamin1},
\cite{bona1}, \cite{grillakis1} and \cite{W1}. Moreover, Kawahara
equation is a model for small-amplitude gravity-capillary waves on
water of a finite depth when the Weber number is close to
$\displaystyle\frac{1}{3}$ (for details, see \cite{schneider}). In
this case, we have a break down when the Weber number is close to
$\displaystyle\frac{1}{3}$. If the Weber number is larger than
$\displaystyle\frac{1}{3}$, this equation has solitary waves just as
the KdV approximation (see
\cite{calvo}).\\
\indent Regarding the stability of solitary waves solutions, we can
mention some contributors. In fact, Angulo in \cite{angulo5} showed
the instability of solitary traveling-wave solutions associated with
the generalized fifth-order KdV equation of the form
\begin{equation}
u_t+u_{xxxxx}+bu_{xxx}=(G(u,u_x,u_{xx}))_x, \label{g5KdV}
\end{equation}
where $ G(q,r,s)=F_q(q,r)-rF_{qr}(q,r)-sF_{rr}(q,r)$ for some
$F(q,r)$ which is homogeneous of degree $p+1$ for some $p>1$, but
the solitary wave was obtained by solving a constrained minimization
problem in $H^2(\mathbb{R})$ which is based on results obtained by
Levandosky (see \cite{levandosky}). The instability of this class of
solitary-wave solutions is determined for $b\neq0$, and it is
obtained by making use of the variational characterization of the
solitary waves and a modification of the theories of instability
established by Shatah\&Strauss \cite{shatah},
Bona\&Souganidis\&Strauss \cite{bona1} and Gon\c calves Ribeiro
\cite{g1}. Levandosky's method in \cite{levandosky}, was also used
by Bridges$\&$Derks \cite{bridges} to show a result of the linear
instability of solitary waves associated with the equation
$(\ref{g5KdV})$. However,
the authors make use of a geometric approach.\\
\indent We recall, from the results of Albert in \cite{albert1}, the
solitary wave
\begin{equation}
u(x,t)=\varphi(x-c_0t)=\displaystyle\mbox{sech}^4\left(x-\frac{12}{35}t\right)
\label{solalbert}\end{equation} where
$c_0\displaystyle=\frac{12}{35}$, is a stable solution of the
Kawahara equation,
$$u_t+uu_x+\displaystyle\frac{13}{420}u_{xxx}-\frac{1}{1680}u_{xxxxx}=0.$$
In this result, the author used the nontrivial polynomials of
Gegenbauer to determine the sign (strictly negative) of the quantity
$I=(\chi,\varphi)_{L^2(\mathbb{R})}$. Here, $\chi\in
L^2(\mathbb{R})$ is such that $\mathcal{L}\chi=\varphi$ (see Theorem
$\ref{teo1}$
in Section $3$).\\
\indent Now, for more general dispersive evolution equations of the
general form
\begin{equation}\label{equa}
u_t+u^pu_x-Mu_x=0,
\end{equation}
an important study of sufficient conditions for the stability was
established by Albert in \cite{albert1} (see also \cite{albert2})
about solitary traveling waves of the form
$
u(x,t)=\varphi(x-ct),
$
for the equation
\begin{equation}\label{trav}
(\mathcal{M} +c)\varphi-\frac{1}{p+1}\varphi^{p+1}=0.
\end{equation}
In (\ref{equa}) (and consequently in $(\ref{trav})$), $p\geq 1$ is
an integer and $\mathcal{M}$ is a Fourier multiplier operator
defined by
\begin{equation}\label{opera1}
\widehat{\mathcal{M}g}(k)=\delta(k)\widehat{g}(k),\;\;\;k\in
\mathbb{R},
\end{equation}
where the symbol $\delta$ is a measurable, locally bounded, even
function on $\mathbb R$ and satisfies that $A_1|k|^{\nu}\leqq
\delta(k)\leqq A_2 (1+|k|)^{\mu} $ for $\nu\leqq \mu$, $ |k|\geq
k_0$, $\delta(k)>b$ for all $k\in \mathbb{R}$ and $A_i>0$. In
\cite{albert1} sufficient conditions were determined to obtain that
the linear, closed, unbounded, self-adjoint operator $\mathcal L:
D(\mathcal L )\to L^2(\mathbb{R})$, defined on a dense subspace of
$L^2(\mathbb{R})$ by
\begin{equation}\label{opera}
\mathcal L \zeta= (\mathcal{M}+c)\zeta-\varphi^p\zeta
\end{equation}
where $\mathcal{M}+c$ is a positive operator, it will have exactly
one negative eigenvalue which is simple and zero is simple with
eigenfunction $\displaystyle\frac{d}{dx}\varphi$. These specific
spectral properties of $\mathcal L$ were obtained provided $\varphi$
is a positive solitary wave satisfying that $\widehat{\varphi}>0$
and $\widehat{\varphi^p}\in PF(2)$ class
defined by Karlin in \cite{karlin1}. \\
\indent In this work, we will show two new explicit families of
stable solitary traveling-wave solutions for the Kawahara and
modified Kawahara equations $(\ref{equakawa})$ and
$(\ref{equamkawa})$ respectively. Such solitary waves are given, in
the case of the Kawahara equation, by
\begin{equation}
\varphi_{\omega}(\xi)=\beta_1\mbox{sech}^2(b\xi)+\lambda_1\mbox{sech}^4(b\xi),
\label{solkawa1}
\end{equation}
where ${\omega}>0$ is the wave-speed and $\beta_1,\lambda_1,b>0$
depending smoothly of $\omega$. In the case of the modified
Kawahara, we have
\begin{equation}
\phi_c(\xi)=\beta_2\mbox{sech}^2(\alpha\xi) \label{solkawa2}
\end{equation}
where $c>0$ is the wave speed with $\alpha$ and $\beta_2>0$ are
parameters which depends smoothly of the wave-speed $\omega$.
However, in this specific case, we cannot obtain the nontrivial
solitary traveling-wave solution associated with the modified
Korteweg-de Vries equation $(\ref{mkdv})$ as $\gamma_2\rightarrow0$,
namely the solitary traveling-wave solution
$g_{\omega}(x)=\displaystyle3\omega\mbox{sech}\left(\frac{\sqrt{\omega}
x}{2}\right)$ associated with the equation $(\ref{mkdv})$. Note that
in $(\ref{solkawa1})$, if $\lambda_1\rightarrow0$ then, we could
expect a profile solitary wave associated
with the KdV equation $(\ref{kdv})$.\\
\indent For both cases, we will use the following conditions that
imply the stability (see \cite{benjamin1}, \cite{bona1},
\cite{grillakis1}, and \cite{W1}):
\begin{equation}\label{propr}
\begin{array}{llll}
(P_0)\;\; \mbox{\rm there is a non-trivial smooth curve of
solutions for (\ref{trav}) of the form},\\
 \;\;\;\;\ \ \ \ c\in I\subseteq \Bbb R\to\varphi_c\in H^{2}(\mathbb{R});\\
(P_1)\;\; \mathcal L\;\;\mbox{\rm has a unique negative
eigenvalue $\lambda$, and which is simple};\\
(P_2)\;\; \mbox{\rm the eigenvalue}\;\; 0 \;\;\mbox{\rm is simple};\\
(P_3)\;\; \frac{d}{dc} \int_{\mathbb{R}} \varphi^2_c(x) dx>0.
\end{array}
\end{equation}
 \indent Therefore, by using conditions in $(\ref{propr})$ we are
 capable to investigate the nonlinear stability of the traveling-wave solutions
of the forms $(\ref{solkawa1})$ and $(\ref{solkawa2})$ for the
Kawahara $(\ref{equakawa})$ and modified Kawahara
$(\ref{equamkawa})$ equations, by using the theory developed by
Albert \cite{albert1} (see also Albert \textit{et al.}
\cite{albert2}). Our stability result is derived from the ideas of
Benjamin\&Bona\&Weinstein\&Grillakis\&Shatah\&Strauss (see
\cite{benjamin1}, \cite{bona1}, \cite{grillakis1} and \cite{W1}).\\
\indent In order to show the current findings, the paper is
organized as follows. Section 2 establishes the notation used in the
body of the paper and well-posedness results for the Kawahara and
modified Kawahara equations. In Section 3 we present the general
theory of stability and the main facts about the paper written by
Albert \cite{albert1}. Section 4 will show the existence of a branch
of solitary traveling waves for the Kawahara equation and the
respective proof of the stability. In Section 5 a branch of solitary
traveling waves will be presented for the equation
$(\ref{equamkawa})$ and the respective proof of the stability for
this case.

\section{Preliminaries and Well-Posedness Results.}
We denote by $\widehat{f}$ the Fourier transform of $f$ in
$\mathbb{R}$, which is defined as
$\widehat{f}(\xi)=\displaystyle\int_{-\infty}^{\infty}\;f(x)e^{-i\xi
x}\,dx.$ Symbol $|f|_{L^p}$ denotes the $L^p(\mathbb{R})$ norm of
$f$, $1\leq p\leq \infty$. In particular, $|\cdot|_{L^2}=\|\cdot\|$
and $|\cdot|_{L^{\infty}}=|\cdot|_\infty$. The inner product of two
elements $f,g\in L^2(\mathbb{R})$ will be denoted by $\langle
f,g\rangle$. We denote by $H^s(\mathbb{R})$, $s\in \mathbb{R}$, the
Sobolev space of all $f$ (tempered distributions) for which the norm
$
\|f\|_{H^s}^2=\displaystyle\int_{-\infty}^{\infty}\;(1+|\xi|^2)^{s}|\widehat
{f}(\xi)|^2\,d\xi$ is finite.\\
\subsection{Well-Posedness Results.}
An interesting result of well-posedness of the Kawahara equation in
$L^2(\mathbb{R})$ is given by Cui$\&$Deng$\&$Tao in \cite{cui}. In
the case of the modified Kawahara the result of well-posedness is
given by Cui\&Tao \cite{cui1}. For both cases, the authors make use
of the techniques of the Bourgain's spaces. These results can be
summarized by the followings Theorems,

\begin{teo} Let $s\geq0$. For each $u_0\in H^s(\mathbb{R})$ there is a
$T>0$ and a unique solution $u\in C([0,T];H^s(\mathbb{R}))$ of the
Kawahara equation $(\ref{equakawa})$. Moreover, the correspondence
$u_0\mapsto u$ is a continuous function between the adequate spaces
\label{wpkawa}\end{teo} \textbf{Proof:} See \cite{cui}.
\begin{flushright}
${\square}$
\end{flushright}

\begin{teo} Let $s\geq2$. For each $u_0\in H^s(\mathbb{R})$ there is a
$T>0$ and a unique solution $u\in C([0,T];H^s(\mathbb{R}))$ of the
modified Kawahara equation $(\ref{equamkawa})$. Moreover, the
correspondence $u_0\mapsto u$ is a continuous function between the
adequate spaces \label{wpmkawa}\end{teo} \textbf{Proof:} See
\cite{cui1}.
\begin{flushright}
${\square}$
\end{flushright}

\section{Stability Theorem and Positivity Properties.} We start
with our definition of stability
\begin{defi}
Let $\varphi$ be a solitary traveling-wave solution of the equation
$(\ref{equakawa})$ (respectively $\ref{equamkawa}$) and consider
$\tau_{r}\varphi(x)=\varphi(x+r)$, $x\in \mathbb{R} $ and
$r\in\mathbb{R}$. We define the set $\Omega_{\varphi}\subset
H^2(\mathbb{R})$, called the orbit generated by $\varphi$, as
$$
\Omega_{\varphi}=\{g;\ g=\tau_r\varphi,\ \mbox{for some}\
r\in\mathbb{R}\}.
$$
And for any $\eta>0$, define the set $U_{\eta}\subset
H^2(\mathbb{R})$ by
$$
U_{\eta}=\left\{f;\
\displaystyle\inf_{g\in\Omega_{\varphi}}||f-g||_{H^2}<\eta\right\}.
$$

\indent With this terminology, we say that $\varphi$ is (orbitally)
stable in $H^2(\mathbb{R})$ by the flow generated by equation
$(\ref{equakawa})$
(respectively $(\ref{equamkawa})$) if,\\
{\rm (i)}  the initial value problem associated with
(\ref{equakawa}) (respectively $(\ref{equamkawa})$) is globally
well-posed in $H^2(\mathbb{R})$ (see Theorems
$\ref{wpkawa}$ and $\ref{wpmkawa}$).\\
{\rm (ii)} For every $\varepsilon>0$, there is $\delta>0$ such that
for all $u_0\in U_{\delta}$, the solution $u$ of $(\ref{equakawa})$
(respectively $(\ref{equamkawa})$) with $u(0,x)=u_0(x)$ satisfies
$u(t)\in U_{\varepsilon}$ for all $t>0$. \label{defi1}
\end{defi}

The proof of the following general stability Theorem can be obtained
by using the techniques given by Benjamin \cite{benjamin1}, Bona
\cite{bona1}, Weinstein \cite{W1} and Grillakis {\it et al.}
\cite{grillakis1}.

\begin{teo}
Let $\varphi$ be a solitary traveling-wave solution of
$(\ref{trav})$ and suppose that part {\rm(i)} of the definition of
stability holds. Suppose also that the operator proceeding of the
equation $(\ref{trav})$,
\begin{equation}
\mathcal{L}\zeta=(\mathcal{M}+c)\zeta-\varphi^p\zeta,
\label{operator}\end{equation} determines that $\mathcal{L}$ has
exactly a unique negative eigenvalue which is simple and zero is a
simple eigenvalue with eigenfunction
$\displaystyle\frac{d}{dx}\varphi_{c}$. Choose $\chi\in
L^2(\mathbb{R})$ such that $\mathcal{L} \chi=\varphi$ and define
$I=(\chi,\varphi)_2$. If $I<0$, then $\varphi$ is stable.
\label{teo1}
\end{teo}
\begin{flushright}
${\square}$
\end{flushright}
\begin{obs}
{\rm (i)} If condition $(P_0)$ in $(\ref{propr})$ holds, we have in
our case that function $\chi$ will be defined as
$\displaystyle\chi=-\frac{d}{d{\omega}}\varphi_{\omega}$ or
$\displaystyle\chi=-\frac{d}{d{c}}\phi_{c}$. Then, it is necessary
to verify that
$\displaystyle\frac{d}{d{\omega}}\|\varphi_{\omega}\|^2>0$ or
$\displaystyle\frac{d}{d{c}}\|\phi_{c}\|^2>0$.\\
{\rm (ii)} The existence of eigenvalues (and as consequence,
eigenfunctions) for the operator $(\ref{operator})$ is guaranteed
from the results contained in \cite{albert3}. \label{obs1}\end{obs}

\indent The main result of the paper in Albert \cite{albert1} (see
also \cite{albert2}) will be presented as follows. Before this, we
need a preliminary definition
\begin{defi}
We say that a function $g:\mathbb{R}\rightarrow\mathbb{R}$ is in the
class $PF(2)$ if \\
{\rm i)} $g(x)>0$, $\forall\ x\in\mathbb{R}$,\\
{\rm ii)}\;$g(x_1-x_2)g(x_2-x_2)-g(x_1-x_2)g(x_2-x_1)>0\ \
\mbox{for}\ x_1<x_2\ \mbox{and}\ x_1<x_2.$ \label{defi2}
\end{defi}

\begin{teo}
Let $\varphi$ be an even  positive solution of (\ref{trav}). Suppose
that $\widehat{\varphi}>0$ and $\mathcal{K}=\widehat{\varphi^p}\in
PF(2)$ discrete, then $\mathcal{L}$ in (\ref{opera}) has exactly a
unique negative eigenvalue which is simple and zero is a simple
eigenvalue with eigenfunction $\displaystyle\frac{d}{dx}\varphi$.
\label{teoprinc}
\end{teo}
\begin{flushright}
${\square}$
\end{flushright}

\section{Existence and Stability of Solitary Traveling-Wave Solutions for the
Kawahara Equation.} In this section we are interested in applying
the theory developed by Albert in \cite{albert1} to obtain the
stability of a specific branch of positive solitary traveling waves
associated with the Kawahara equation whose statements was presented
in the previous section.

\subsection{Existence of Solitary Traveling-Wave Solutions.}
In this subsection we establish the existence of solitary
traveling-wave solutions related to the Kawahara equation given by,

\begin{equation}
u_t+uu_x+u_{xxx}-\gamma_1 u_{xxxx}=0. \label{kawa1}
\end{equation}
In fact, let $u(x,t)=\varphi_{\omega}(x-\omega t)$ be a solitary
traveling-wave solution associated with $(\ref{kawa1})$.
Substituting this form in the equation $(\ref{kawa1})$ we obtain,
after integration, that
\begin{equation}
\displaystyle-{\omega}\varphi_{\omega}+\frac{1}{2}\varphi_{\omega}^2
+\varphi_{\omega}''-\gamma_1\varphi_{\omega}''''=0,
\label{solitkawa}
\end{equation}
where ${\omega}\in\mathbb{R}$.\\
\indent Next, we consider
\begin{equation}\varphi_{\omega}(\xi)=\beta_1\mbox{sech}^2(b\xi)+\lambda_1\mbox{sech}^4(b\xi),
\label{solusolit}\end{equation} where $\beta_1,\lambda_1$ and $b>0$
a smooth solution for $(\ref{solitkawa})$. By using \textit{Maple
program}, the following nonlinear system is obtained,
\begin{equation}
\left\{\begin{array}{lll}
\lambda_1-1680b^4\gamma_1=0\\\\
\displaystyle-\frac{1}{2}b^2+\gamma_1 b^4+\frac{1}{16}\omega=0\\\\
\displaystyle 240\gamma_1 b^4\beta_1-512\gamma_1
b^4\lambda_1-32\lambda_1b^2-12\beta_1 b^2
+\beta_1^2-2\omega\lambda_1=0\\\\
\displaystyle 2080\gamma_1 b^4\lambda_1-240\gamma_1\beta_1
b^4-40\lambda_1b^2+2\beta_1\lambda_1=0.\end{array}\right.\label{sistsolkawa}
\end{equation}

After some calculations, $(\ref{sistsolkawa})$ boils down in a
simple system as

\begin{equation}
\left\{\begin{array}{lll}
\displaystyle b^2-\frac{\lambda_1}{840}-\frac{\omega}{8}=0\\\\
26\lambda_1+39\beta_1-840b^2=0\\\\
3\lambda_1\beta_1-32\lambda_1^2-672\lambda_1b^2-252\beta_1b^2+21\beta_1^2
-42\omega\lambda_1=0
\end{array}\right.\label{sistsolkawa1}
\end{equation}

System $(\ref{sistsolkawa1})$ can be dropped in terms of $\lambda_1$
and $\omega$ as

\begin{equation}\label{lambda-omega}
\displaystyle-\frac{2023210}{169}\omega\lambda_1-\frac{862463}{507}\lambda_1^2
+\frac{797475}{169}\omega^2=0.
\end{equation}
Then, we discover $\lambda_1$ in terms of $\omega$
\begin{equation}\label{lambda}
\lambda_1(\omega)=\displaystyle
105\left(-\frac{4129}{123209}+\frac{546}{123209}\sqrt{70}\right)\omega.
\end{equation}
where we can conclude that $\lambda_1(\omega)>0$ and
$\lambda_1'(\omega)>0$ for all $\omega>0$.\\
\indent Further, from $(\ref{sistsolkawa1})$ and $(\ref{lambda})$ we
have

\begin{equation}\label{beta}\begin{array}{lll}
\beta_1(\omega)&=&\displaystyle\frac{105\omega}{39}-\frac{25\lambda_1}{39}\\\\
&=&\displaystyle\frac{\left(609630-36750\sqrt{70}\right)\omega}{123209}.
\end{array}
\end{equation}
Therefore, we get $\beta_1(\omega)>0$ and $\beta_1'(\omega)>0$ for
all $\omega>0$.\\
\indent Finally, we can find $b$ in term of $\omega$ as

\begin{equation}\label{b}
b(\omega)=\frac{\sqrt{123209}\sqrt{(59540+273\sqrt{70})\omega}}{246418}
\end{equation}
and we have $b(\omega)>0$ and $b'(\omega)>0$ for all $\omega>0$.

Next, from $(\ref{lambda})$, $(\ref{beta})$ and $(\ref{b})$ we
deduce that
\begin{equation}\label{existsolit}
\omega\in(0,+\infty)\mapsto\varphi_{\omega}\in H^n(\mathbb{R})\ \
\mbox{is smooth}\ \ \mbox{for all}\ \ n\in\mathbb{N}.
\end{equation}

\subsection{Stability of Solitary Traveling-Wave Solutions.}
We have the following Theorem of stability

\begin{teo}
The smooth branch of solutions $\varphi_{\omega}$ obtained in
$(\ref{existsolit})$ is orbitally stable in $H^2(\mathbb{R})$ by the
flow of the Kawahara equation since $\omega>0$.
\label{teoestsolit}\end{teo} \textbf{Proof:} First of all we wish to
determine the behavior of the first two eigenvalues associated with
the operator
$$\mathcal{L}=\displaystyle\gamma_1\frac{d^4}{dx^4}-\frac{d^2}{dx^2}
+\omega-\varphi_{\omega},$$ by utilizing the theory developed by
Albert in \cite{albert1}. In fact, since the kernel
$\mathcal{K}=\widehat{\varphi_{\omega}}$ belongs to the $PF(2)$
continuous case from the Lemma $10$ in \cite{albert2}, it is
necessary to show that the symbol $\delta(z)$ associated with the
linear operator
$\mathcal{M}=\displaystyle\gamma_1\frac{d^4}{dx^4}-\frac{d^2}{dx^2}$
satisfies the properties in $(\ref{opera1})$. Indeed, since
$\delta(z)=\widehat{Mu}(z)=(\gamma_1|z|^4+|z|^2)\widehat{u}(z)$ for
all $z\in\mathbb{R}$, the properties are verified.\\
\indent Next, from Theorem $\ref{teo1}$ and Remark $\ref{obs1}$-{\rm
(i)} we calculate
$\displaystyle\frac{d}{d{\omega}}||\varphi_{\omega}||^2$, where
$\varphi_{\omega}(\xi)=\beta_1\mbox{sech}^2(b\xi)+\lambda_1\mbox{sech}^4(b\xi)$
and $\omega>0$. In fact
$$
\begin{array}{lll}
||\varphi_{\omega}||^2&=&\displaystyle\frac{\beta_1^2}{b}
\int_{\mathbb{R}}\mbox{sech}^4(x)dx+\frac{2\beta_1\lambda_1}{b}
\int_{\mathbb{R}}\mbox{sech}^6(x)dx+\frac{\lambda_1^2}{b}
\int_{\mathbb{R}}\mbox{sech}^8(x)dx\\\\
&=&\displaystyle\frac{4\beta_1^2}{3b}+\frac{32\beta_1
\lambda_1}{15b}+\frac{32\lambda_1^2}{35b}
\end{array}
$$

Since $\displaystyle\frac{\beta_1^2}{b}=M_1\omega^{3/2},\ \
\frac{\beta_1\lambda_1}{b}=M_2\omega^{3/2},\ \
\frac{\lambda_1^2}{b}=M_3\omega^{3/2},$ where $M_i$, $i=1,2,3$ are
positive constants obtained from $(\ref{lambda})$, $(\ref{beta})$
and $(\ref{b})$ we deduce that $d''(\omega)>0$, for all $\omega>0$.
\begin{flushright}
${\square}$
\end{flushright}
\section{Existence and Stability of Solitary Traveling-Wave
Solutions for the modified Kawahara Equation.} This section is
concerned to proof the existence and stability of solitary
traveling-wave solutions for the modified Kawahara equation
\begin{equation}
\displaystyle-c\phi_c+\phi_c^3+\phi_c''-\gamma_2\phi_c''''=0,
\label{travmkawa}\end{equation} where $\gamma_2>0$.

To prove the existence, let
$\phi_c(\xi)=\beta_2\mbox{sech}^2(\alpha\xi)$ be a solitary
traveling-wave solution for the equation $(\ref{travmkawa})$. If we
substitute this $\phi_c$ into $(\ref{travmkawa})$ we obtain after
some calculations $\beta_2=6\alpha$ and
$\alpha=\displaystyle\frac{\sqrt{5c}}{4}$. Therefore for all $c>0$
we have,

\begin{equation}
\phi_c(\xi)=\displaystyle\frac{3\sqrt{5c}}{2}\mbox{sech}^2\left(\frac{\sqrt{5c}}{4}\xi\right).
\label{solmkawa}\end{equation}

\subsection{Stability of Solitary Traveling-Wave Solutions.}
In this subsection we are interested in applying the theory in
Section 3 to obtain the stability of the smooth branch of positive
solitary traveling waves obtained in the last subsection, associated
to the modified Kawahara equation. In fact, our intention can be
summarized in the following theorem,

\begin{teo}
The branch of solutions $\phi_{c}$ given by $(\ref{solmkawa})$ is
orbitally stable in $H^2(\mathbb{R})$ by the flow of the modified
Kawahara equation for all $c>0$.\label{teoestsolit2}\end{teo}
\textbf{Proof:} First of all, we note clearly that
$\widehat{\phi_c}>0$. Lemma 10 in \cite{albert2}, also shows that
$\widehat{\phi_c^2}$ belongs to the $PF(2)$ class and therefore the
properties $(P_1)$ and $(P_2)$ in $(\ref{propr})$ are satisfied for
the operator
$\mathcal{L}=\displaystyle\gamma_2\frac{d^4}{dx^4}-\frac{d^2}{dx^2}+c-3\phi_c^2,$
associated with the equation $(\ref{travmkawa})$. Moreover, the
linear operator
$\mathcal{M}=\displaystyle\gamma_2\frac{d^4}{dx^4}-\frac{d^2}{dx^2}$,
satisfies the properties required in $(\ref{opera1})$ by the same
arguments seen in the proof of the Theorem $\ref{teoestsolit}$. It
remains for us to calculate the quantity
$\displaystyle\frac{d}{dc}||\phi_c||^2$. In fact, since
$$\displaystyle\int_{\mathbb{R}}\phi_c(\xi)^2d\xi=\frac{45\sqrt{c}}{\sqrt{5}}\int_{\mathbb{R}}\mbox{sech}^4(y)dy
=\frac{60\sqrt{c}}{\sqrt{5}},$$ we have
$\displaystyle\frac{d}{dc}||\phi_c||^2>0$. This argument shows the
theorem.
\begin{flushright}
${\square}$
\end{flushright}

\end{document}